\documentclass{amsart}

\usepackage{graphicx,pstricks,pst-plot}

\newtheorem{theorem}{Theorem}[section]
\newtheorem{lemma}[theorem]{Lemma}
\newtheorem{proposition}[theorem]{Proposition}
\newtheorem{corollary}[theorem]{Corollary}

\theoremstyle{definition}

\theoremstyle{remark}
\newtheorem{remark}[theorem]{Remark}

\numberwithin{equation}{section}

\title[A Generalized Maximum Principle for Yau's Operator]{A Generalized Maximum Principle for \\Yau's Square Operator, with Applications \\to the Steady State Space}

\author{A. Caminha}
\address{Departamento de Matem\'atica, Universidade Federal do Cear\'a,
Fortaleza, Cear\'a, Brazil. 60455-760}
\email{antonio.caminha@gmail.com}
\thanks{The authors wish to thank professors L. J. Al\'{\i}as and A. G. Colares for some valuable comments.}

\author{H. F. de Lima}
\address{Departamento de Matem\'atica e Estat\'{\i}stica, Universidade
Federal de Campina Grande, Campina Grande, Para\'{\i}ba, Brazil.
58109-970} \email{henrique@dme.ufcg.edu.br}
\thanks{The second author is partially supported by CAPES}


\begin{document}
\maketitle

\begin{abstract}
We derive, for the square operator of Yau, an analogue of the Omori-Yau maximum principle for the Laplacian.
We then apply it to obtain nonexistence results concerning complete spacelike hypersurfaces
with constant higher order mean curvature in the Steady State space.
\end{abstract}

\section{Introduction}

The interest in the study of spacelike hypersurfaces in Lorentz manifolds ({\em spacetimes}) has increased very much
in recent years, from both the physical and mathematical points of view. A basic question on this topic is the existence
and uniqueness of spacelike hypersurfaces with some reasonable geometric properties, like the vanishing of
the mean curvature, for instance. A first relevant result in this direction was the proof of the Calabi-Bernstein
conjecture for maximal hypersurfaces (that is, hypersurfaces with vanishing mean curvature) in Lorentz-Minkowski
space, given by Cheng an Yau in~\cite{Cheng:76}. As for the case of de Sitter space, Goddard in~\cite{Goddard:77}
conjectured that every complete spacelike hypersurface with constant mean curvature in de Sitter space should be
totally umbilical. Although the conjecture turned out to be false in its original form, it motivated a great deal
of work of several authors trying to find a positive answer to the conjecture under appropriated additional
hypotheses (see, for example, ~\cite{Akutagawa:87} and ~\cite{Montiel1:88}).

More recently, Al\'{\i}as, Brasil and Colares, in~\cite{AliasBrasil:03}, developed general Minkowski-type formulae for
compact spacelike hypersurfaces immersed into conformally stationary spacetimes, that is, spacetimes endowed with a
timelike conformal vector field; then, they applied these formulae to the study of the umbilicity of compact
spacelike hypersurfaces under appropriate conditions on their $r$-mean curvatures. Furthermore, the first author
in~\cite{Caminha:06} computed $L_r(S_r)$ for a spacelike hypersurface $\Sigma^n$ immersed in a spacetime
$\overline M^{n+1}$ of constant sectional curvature, applying the resulting formula to study both
$r$-maximal spacelike hypersurfaces of $\overline M$, and, in the presence of a constant higher order mean
curvature, constraints on the sectional curvature of $\Sigma$ that also suffice to guarantee the umbilicity of it.
Here, by $L_r$ we mean the linearization of the second order differential operator associated to the $r$-th elementary
symmetric function $S_r$ on the eigenvalues of the second fundamental form of such immersion (cf. section
\ref{sec:preliminaries}). Let us also remark that Al\'{\i}as and Colares in ~\cite{AliasColares:06} studied the problem
of uniqueness for spacelike hypersurfaces with constant higher order mean curvature in {\em Generalized Robertson-Walker}
spacetimes. Their approach is based on the use of the Newton transformations $P_r$ (and their associated
differential operators $L_r$) and the abovementioned Minkowski formulae for spacelike hypersurfaces.

Returning to the case of complete noncompact spacelike hypersurfaces, the first author in~\cite{Caminha:07} used the
standard formula for the Laplacian of the squared norm of the second fundamental form and the Omori-Yau maximum
principle to classify complete spacelike hypersurfaces with constant mean curvature in a spacetime of nonnegative
constant sectional curvature, under appropriate bounds on the scalar curvature. For the de Sitter space, Brasil Jr.,
Colares and Palmas also used the Omori-Yau maximum principle in ~\cite{Brasil:03} to characterize the hyperbolic
cylinders as the only complete hypersurfaces in the de Sitter space with constant mean curvature, nonnegative Ricci
curvature and having at least two ends (see also ~\cite{Brasil:01} for the case of the scalar curvature).

The discussion of related questions involving higher order mean curvatures faces a first difficulty: there is
no corresponding version of maximum principle for the appropriate second order partial differential operators. Therefore,
we begin this paper by overcoming this obstacle. More precisely, if $\Phi$ is a field of self-adjoint linear maps
on a spacelike hypersurface $\Sigma^n$ of a spacetime $\overline{M}^{n+1}$, and $f\in\mathcal D(\Sigma)$,
we consider the square operator
$$\square f={\rm tr}(\Phi\,\text{Hess}\,f)$$
of Yau~\cite{Cheng:77}; when
$\Sigma$ is complete, and under certain conditions on $\Phi$, we develop for $\square$ an analogue of the Omori-Yau
maximum principle for the Laplacian (see Corollary~\ref{corol:functions bounded from below}).
With the aid of a suitable corollary of it (cf. Proposition~\ref{prop:Akutagawa for square operator}), we obtain
nonexistence results on complete noncompact spacelike hypersurfaces $\Sigma^n$ into the half $\mathcal H^{n+1}$ of
the de Sitter space having one constant higher order mean curvature (cf. theorems~\ref{thm:application 1},
~\ref{thm:application 2} and~\ref{thm:application 3}).

We observe that $\mathcal H^{n+1}$, the so-called {\em Steady State space}, appears naturally in physical context
as an exact solution for the Einstein equations, being a cosmological model where matter is supposed to travel along
geodesics normal to horizontal hyperplanes (slices); these, in turn, serve as the initial data for the Cauchy problem
associated to those equations (cf.~\cite{Hawking:73}, chapter $5$).

This paper is organized in the following manner: in section~\ref{sec:preliminaries} we set
notation and recall a few results which will be needed later; section~\ref{sec:The generalized maximum principle}
is devoted to the statement and proof of the maximum principle and its corollaries; applications are collected in
section~\ref{sec:applications}.

\section{Preliminaries}\label{sec:preliminaries}

Let $\overline M^{n+1}$ be a connected Lorentz manifold with metric $\overline g=\langle\,\,,\,\,\rangle$, and
Levi-Civitta connection $\overline\nabla$. We recall (cf.~\cite{O'Neill:83}) that a vector field
$X\in\mathcal X(\overline M)$ is said to be
{\em timelike} if $\langle X,X\rangle<0$ on $\overline M$; {\em spacelike} if $\langle X,X\rangle>0$ on $\overline M$;
a {\em unit} vector field if $\langle X,X\rangle=\pm 1$ on $\overline M$.

A vector field $V$ on $\overline M^{n+1}$ is said to be {\em conformal} if
\begin{equation}
\mathcal L_V\langle\,\,,\,\,\rangle=2\phi\langle\,\,,\,\,\rangle
\end{equation}
for some function $\phi\in C^{\infty}(\overline M)$, where $\mathcal L$ stands for the Lie
derivative of the metric of $\overline M$. The function $\phi$ is called the
{\em conformal factor} of $V$.

Since $\mathcal L_V(X)=[V,X]$ for all $X\in\mathcal X(\overline M)$, it follows from the
tensorial character of $\mathcal L_V$ that $V\in\mathcal X(\overline M)$ is conformal if
and only if
\begin{equation}\label{eq:1.1}
\langle\overline\nabla_XV,Y \rangle+\langle X,\overline\nabla_YV\rangle=2\phi\langle X,Y\rangle,
\end{equation}
for all $X,Y\in\mathcal X(\overline M)$. In particular, $V$ is a Killing vector field
relatively to $\overline g$ if and only if $\phi\equiv 0$.

In all that follows, we consider spacelike immersions $\psi:\Sigma^n\rightarrow\overline M^{n+1}$, namely,
isometric immersions from a connected, $n-$dimensional orientable Riemannian manifold $\Sigma$ into $\overline M$.
We let $\nabla$ denote the Levi-Civita connection of $\Sigma$.

Let us orient $\Sigma$ by the choice of a unit normal vector field $N$ on it, and $A$ denote the corresponding shape
operator. At each $p\in\Sigma$, $A$ restricts to a self-adjoint linear map $A_p:T_p\Sigma\rightarrow T_p\Sigma$.
For $1\leq r\leq n$, let $S_r(p)$ denote the $r-$th elementary symmetric function on the eigenvalues of $A_p$;
this way one gets $n$ smooth functions $S_r:\Sigma^n\rightarrow\mathbb R$, such that
$$\det(tI-A)=\sum_{k=0}^n(-1)^kS_kt^{n-k},$$
where $S_0=1$ by definition. If $p\in\Sigma$ and $\{e_k\}$ is a basis of $T_p\Sigma$ formed by
eigenvectors of $A_p$, with corresponding eigenvalues $\{\lambda_k\}$, one immediately sees
that
$$S_r=\sigma_r(\lambda_1,\ldots,\lambda_n),$$
where $\sigma_r\in\mathbb R[X_1,\ldots,X_n]$ is the $r-$th elementary symmetric polynomial
on the indeterminates $X_1,\ldots,X_n$. In particular, if $|A|^2$ stands for ${\rm tr}(A^2)$
then it is immediate to check that
\begin{equation}\label{eq:S1 S2 and A2}
2S_2+|A|^2=S_1^2.
\end{equation}

For $1\leq r\leq n$, one defines the $r-th$ mean curvature $H_r$ of $\psi$ by
$$H_r=\frac{(-1)^r}{{n\choose r}}S_r=\frac{1}{{n\choose r}}\sigma_r(-\lambda_1,\ldots,-\lambda_n).$$
In particular, $H_1=H$ is the mean curvature of $x$. It is a classical fact that such functions satisfy
a very useful set of inequalities, usually reffered to as Newton's inequalities (see~\cite{Hardy:89}).
It turns out, however, that such inequalities remain true for arbitrary real numbers.
For future reference, we collect them here. A proof can be found in~\cite{Caminha:06}, proposition 1.

\begin{proposition}\label{prop:Newton's inequalities}
Let $n>1$ be an integer, and $\lambda_1,\ldots,\lambda_n$ be real numbers.
Define, for $0\leq r\leq n$, $S_r=S_r(\lambda_i)$ as above, and
$H_r=H_r(\lambda_i)={n\choose r}^{-1}S_r(\lambda_i)$.
\begin{enumerate}
\item[(a)] For $1\leq r<n$, one has $H_r^2\geq H_{r-1}H_{r+1}$. Moreover, if
equality happens for $r=1$ or for some $1<r<n$, with $H_{r+1}\neq 0$ in this case,
then $\lambda_1=\cdots=\lambda_n$.
\item[(b)] If $H_1,H_2,\ldots,H_r>0$ for some $1<r\leq n$, then
$H_1\geq\sqrt{H_2}\geq\sqrt[3]{H_3}\geq\cdots\geq\sqrt[r]{H_r}$. Moreover, if equality
happens for some $1\leq j<r$, then $\lambda_1=\cdots=\lambda_n$.
\item[(c)] If, for some $1\leq r<n$, one has $H_r=H_{r+1}=0$, then
$H_j=0$ for all $r\leq j\leq n$. In particular, at most $r-1$ of the $\lambda_i$ are
different from zero.
\end{enumerate}
\end{proposition}

When the ambient space $\overline M$ has constant sectional curvature $c$, Gauss
equation allows one to immediately check that the scalar curvature $R$ of $\Sigma$ relates to
$H_2$ in the following manner:
\begin{equation}\label{eq:S2 and scalar curv}
R=n(n-1)(c-H_2).
\end{equation}

For $0\leq r\leq n$ one defines the $r-$th Newton transformation $P_r$ on $\Sigma$ by
setting $P_0=I$ (the identity operator) and, for $1\leq r\leq n$, via the recurrence relation
\begin{equation}\label{eq:Newton operators}
P_r=(-1)^rS_rI+AP_{r-1}.
\end{equation}
A trivial induction shows that
$$P_r=(-1)^r(S_rI-S_{r-1}A+S_{r-2}A^2-\cdots+(-1)^rA^r),$$
so that Cayley-Hamilton theorem gives $P_n=0$. Moreover, since $P_r$ is a polynomial in
$A$ for every $r$, it is also self-adjoint and commutes with $A$. Therefore, all bases
of $T_p\Sigma$ diagonalizing $A$ at $p\in\Sigma$ also diagonalize all of the $P_r$ at $p$.
Let $\{e_k\}$ be such a basis. Denoting by $A_i$ the restriction of $A$ to
$\langle e_i\rangle^{\bot}\subset T_p\Sigma$, it is easy to see that
$$\det(tI-A_i)=\sum_{k=0}^{n-1}(-1)^kS_k(A_i)t^{n-1-k},$$
where
$$S_k(A_i)=\sum_{\stackrel{1\leq j_1<\ldots<j_k\leq n}{j_1,\ldots,j_k\neq i}}\lambda_{j_1}\cdots\lambda_{j_k}.$$

With the above notations, it is also immediate to check that $P_re_i=(-1)^rS_r(A_i)e_i$,
and hence (lemma 2.1 of~\cite{Barbosa:97})
\begin{enumerate}
\item[(a)] $S_r(A_i)=S_r-\lambda_iS_{r-1}(A_i)$;
\item[(b)] ${\rm tr}(P_r)=(-1)^r\sum_{i=1}^nS_r(A_i)=(-1)^r(n-r)S_r=b_rH_r$;
\item[(c)] ${\rm tr}(AP_r)=(-1)^r\sum_{i=1}^n\lambda_iS_r(A_i)=(-1)^r(r+1)S_{r+1}=-b_rH_{r+1}$;
\item[(d)] ${\rm tr}(A^2P_r)=(-1)^r\sum_{i=1}^n\lambda_i^2S_r(A_i)=(-1)^r(S_1S_{r+1}-(r+2)S_{r+2})$,
\end{enumerate}
where $b_r=(n-r){n\choose r}$.

The next two results will be extremely useful in section~\ref{sec:applications}.

\begin{proposition}[proposition 1.5 of~\cite{Hounie:99}]\label{prop:ellipiticy_for_null_Sr}
With respect to a spacelike immersion $\psi:\Sigma^n\rightarrow\overline M^{n+1}$,
\begin{enumerate}
\item[(a)] if $H_r=0$ on $\Sigma$, then $P_{r-1}$ is semi-definite on $\Sigma$.
\item[(b)] if $H_r=0$ and $H_{r+1}\neq 0$ on $\Sigma$, then $P_{r-1}$ is definite on $\Sigma$.
\end{enumerate}
\end{proposition}

If $p\in\Sigma$ is such that all eigenvalues of $A_p$ are negative, we say that $p$ is an
{\em elliptic point} of $\Sigma$.

\begin{proposition}[proposition 3.2 of~\cite{Barbosa:97}]\label{prop:ellipiticy_for_positive_Sr}
With respect to a spacelike immersion $\psi:\Sigma^n\rightarrow\overline M^{n+1}$,
if $H_r>0$ on $\Sigma$ and $\psi$ has an elliptic point, then $P_{r-1}$ is positive definite on $\Sigma$.
\end{proposition}

Associated to each Newton transformation $P_r$ one has the second order linear differential
operator $L_r:\mathcal D(\Sigma)\rightarrow\mathcal D(\Sigma)$, given by
$$L_r(f)={\rm tr}(P_r\,\text{Hess}\,f).$$
Therefore, for $f,g\in\mathcal D(\Sigma)$,
it follows from the properties of the Hessian of functions that
\begin{equation}\label{eq:Lr_of_a_product}
L_r(fg)=fL_r(g)+gL_r(f)+2\langle P_r\nabla f,\nabla g\rangle.
\end{equation}

\section{The generalized maximum principle}\label{sec:The generalized maximum principle}

Let $\Sigma^n$ be a complete $n-$dimensional Riemannian manifold.
Let also $\Phi:T\Sigma\rightarrow T\Sigma$ denote a field of self adjoint linear transformations on $\Sigma$.
We consider the second order linear differential operator $\square:\mathcal D(\Sigma)\rightarrow\mathcal D(\Sigma)$
by setting
\begin{equation}
\square f=\,\text{tr}\,(\phi\,\text{Hess}\,f).
\end{equation}

For fixed $p\in\Sigma$, let $\rho(x)=\rho_p(x)=d(x,p)$ be the distance function from $p$ and
$C_m(p)$ denote the cut locus of $p$. Set also
\begin{equation}\label{eq:definition of K}
K(x)=s_c'(\rho(x))s_c(\rho(x)){\rm tr}(\Phi_x),
\end{equation}
where $s_c:[0,+\infty]\rightarrow\mathbb R$ is defined by
$$s_c(t)=\left\{\begin{array}{l}
\frac{\sinh(\sqrt{-c}t)}{\sqrt{-c}t},\,\,\text{if}\,\,c<0\\
t,\,\,\text{if}\,\,c=0\\
\frac{\sin(\sqrt ct)}{\sqrt ct},\,\,\text{if}\,\,c>0
\end{array}\right. .$$

\begin{lemma}\label{lemma:bounds for K with bounded sectional curvature}
If $\Phi$ is positive semi-definite on $\Sigma$ and $\Sigma$ has sectional curvature
$K_{\Sigma}\geq c$ then, for all $x\in\Sigma\setminus C_m(p)$, one has $\square\rho(x)\leq K(x)$.
\end{lemma}

\begin{proof} Let $\gamma:[0,l]\rightarrow\Sigma$ be the only minimizing normalized geodesic joining $p$ to
$x$, with length $l=\rho(x)$. Decompose any unit vector $u\in T_x\Sigma$ as $u=v+w$, where
$u$ is collinear with $\gamma'(l)$ and $w\perp\gamma'(l)$. Then $|v|^2+|w|^2=1$ and, at $x$,
\begin{eqnarray*}
\text{Hess}\,\rho(u,u)&=&\text{Hess}\,\rho(v,v)+2\text{Hess}\,\rho(v,w)+\text{Hess}\,\rho(w,w)\\
&=&\langle\nabla_v\gamma',v\rangle+2\langle\nabla_v\gamma',w\rangle+\text{Hess}\,\rho(w,w)\\
&=&\text{Hess}\,\rho(w,w).
\end{eqnarray*}
It follows from the Hessian comparison theorem and from the the characterization
of Jacobi fields in spaces of constant sectional curvature that if $K_{\Sigma}\geq c$ then, at $x$,
$$\text{Hess}\,\rho(w,w)\leq s_c'(\rho)s_c(\rho)|w|^2\leq s_c'(\rho)s_c(\rho).$$

Now take a moving frame $\{e_1,\ldots,e_n\}$ on a neighborhood of $x$, diagonalizing $\Phi$
at $x$, with $\Phi(e_i)=\lambda_ie_i$. Then, one has at $x$
\begin{eqnarray*}
\square\rho&=&{\rm tr}(\Phi\,\text{Hess}\,\rho)=\sum_i\lambda_i\,\text{Hess}\,\rho(e_i,e_i)\\
&\leq&\sum_i\lambda_is_c'(\rho)s_c(\rho)=s_c'(\rho)s_c(\rho){\rm tr}(\Phi).
\end{eqnarray*}
\end{proof}

\begin{theorem} Let $\Sigma$ be a complete Riemannian manifold with sectional curvature $K_{\Sigma}\geq c$,
and $f\in\mathcal D(\Sigma)$ be a function bounded from above. If $\Phi$ is positive semi-definite
at every $x\in\Sigma$ then, for every $p\in\Sigma$, there exists a sequence $(p_k)_{k\geq 1}$ in $\Sigma$
such that
\begin{equation}\label{eq:bounds for f}
\lim_{k\rightarrow+\infty}f(p_k)=\sup_{\Sigma}f,
\end{equation}

\begin{equation}\label{eq:bounds for the gradient of f}
|\nabla f(p_k)|=\frac{2(f(p_k)-f(p)+1)\rho(p_k)}{k(\rho(p_k)^2+2)\log(\rho(p_k)^2+2)}
\end{equation}
and
\begin{eqnarray}\label{eq:bounds for square of f}
\square f(p_k)&\leq&\frac{4{\rm tr}(\Phi_{p_k})\rho(p_k)^2(f(p_k)-f(p)+1)}
{k^2(\rho(p_k)^2+2)^2\log(\rho(p_k)^2+2)^2}\\
&&+\frac{2(f(p_k)-f(p)+1)}{k(\rho(p_k)^2+2)\log(\rho(p_k)^2+2)}\left\{{\rm tr}(\Phi_{p_k})+\rho(p_k)K(p_k)\right\},\nonumber
\end{eqnarray}
where $K$ is given as in (\ref{eq:definition of K}).
\end{theorem}

\begin{proof} The proof parallels that of the classical Omori-Yau maximum principle in~\cite{Yau:75}.
For positive integer $k$, let
$$g(x)=\frac{f(x)-f(p)+1}{[\log(\rho(x)^2+2)]^{1/k}}.$$
One has that $g$ is continuous, $g(p)=\frac{1}{(\log 2)^{1/k}}>0$ and, since $f$ is bounded
above,
$$\limsup_{\rho(x)\rightarrow+\infty}g(x)\leq 0.$$
Therefore, $g$ attains its maximum at some $p_k\in\Sigma$. In particular, $f(p_k)-f(p)+1>0$.
One now has to consider two cases separately: $p_k\notin C_m(p)$ and $p_k\in C_m(p)$. Here, we treat
only the first case; for the second one and the conclusion of the proof of the theorem, copy the
corresponding steps in~\cite{Yau:75}.

Suppose $p_k\notin C_m(p)$. Since (omitting $x$ for clarity)
\begin{equation}\label{eq:v of g}
v(g)=\frac{v(f)}{[\log(\rho^2+2)]^{1/k}}-\frac{2(f-f(p)+1)\rho v(\rho)}{k(\rho^2+2)[\log(\rho^2+2)]^{1/k+1}},
\end{equation}
one gets at $p_k$
$$0=\nabla g=\frac{\nabla f}{[\log(\rho^2+2)]^{1/k}}-\frac{2(f-f(p)+1)\rho\nabla\rho}{k(\rho^2+2)[\log(\rho^2+2)]^{1/k+1}},$$
from where (\ref{eq:bounds for the gradient of f}) follows.

For the estimate on $\square f$, it follows from (\ref{eq:v of g}) that
\begin{eqnarray*}
v(v(g))&=&\frac{v(v(f))}{[\log(\rho^2+2)]^{1/k}}-\frac{2\rho v(f)v(\rho)}{k(\rho^2+2)[\log(\rho^2+2)]^{1/k+1}}\\
&&-\frac{2\left\{\rho v(f)v(\rho)+(f-f(p)+1)[v(\rho)^2+\rho v(v(\rho))]\right\}}{k(\rho^2+2)[\log(\rho^2+2)]^{1/k+1}}\\
&&+\frac{4(f-f(p)+1)\rho^2v(\rho)^2}{k(\rho^2+2)^2[\log(\rho^2+2)]^{1/k+2}}\left(\frac{1}{k}+1+\log(\rho^2+2)\right).
\end{eqnarray*}

Now take a moving frame $\{e_1,\ldots,e_n\}$ on a neighborhood of $p_k$, geodesic at $p_k$
and diagonalizing $\Phi$ at $p_k$, with $\Phi(e_i)=\lambda_ie_i$. Then, one has at $p_k$
$$\square f=\sum_i\lambda_ie_i(e_i(f)).$$
On the other hand, since $\text{Hess}\,f_{p_k}\leq 0$ and $\Phi_{p_k}\geq 0$, one has
$\square g={\rm tr}(\Phi\,\text{Hess}\,g)\leq 0$ at $p_k$, and it follows at once from the above
computations that
\begin{eqnarray*}
0\geq\square g&=&\frac{\square f}{[\log(\rho^2+2)]^{1/k}}-\frac{4\rho\langle\Phi\nabla f,\nabla\rho\rangle}{k(\rho^2+2)[\log(\rho^2+2)]^{1/k+1}}\\
&&-\frac{2(f-f(p)+1)(\langle\Phi\nabla\rho,\nabla\rho\rangle+\rho\square\rho)}{k(\rho^2+2)[\log(\rho^2+2)]^{1/k+1}}\\
&&+\frac{4(f-f(p)+1)\rho^2\langle\Phi\nabla\rho,\nabla\rho\rangle}{k(\rho^2+2)^2[\log(\rho^2+2)]^{1/k+2}}\left(\frac{1}{k}+1+\log(\rho^2+2)\right).
\end{eqnarray*}
One also has at $p_k$ that
$$\langle\Phi\nabla f,\nabla\rho\rangle=\frac{2(f-f(p)+1)\rho\langle\Phi\nabla\rho,\nabla\rho\rangle}{k(\rho^2+2)\log(\rho^2+2)},$$
from where, substituting into the above and taking into account
lemma~\ref{lemma:bounds for K with bounded sectional curvature}, we get at $p_k$
\begin{eqnarray*}
\square f&\leq&\frac{8(f-f(p)+1)\rho^2\langle\Phi\nabla\rho,\nabla\rho\rangle}{k^2(\rho^2+2)^2[\log(\rho^2+2)]^2}+\frac{2(f-f(p)+1)(\langle\Phi\nabla\rho,\nabla\rho\rangle+\rho K)}{k(\rho^2+2)\log(\rho^2+2)}\\
&&-\frac{4(k+1)(f-f(p)+1)\rho^2\langle\Phi\nabla\rho,\nabla\rho\rangle}{k^2(\rho^2+2)^2[\log(\rho^2+2)]^2}-\frac{4(f-f(p)+1)\rho^2\langle\Phi\nabla\rho,\nabla\rho\rangle}{k(\rho^2+2)^2\log(\rho^2+2)}\\
&=&\frac{2(f-f(p)+1)(\langle\Phi\nabla\rho,\nabla\rho\rangle+\rho K)}{k(\rho^2+2)\log(\rho^2+2)}\\
&&+\frac{4(f-f(p)+1)\rho^2\langle\Phi\nabla\rho,\nabla\rho\rangle}{k^2(\rho^2+2)^2[\log(\rho^2+2)]^2}[2-(k+1)-k\log(\rho^2+2)]\\
&\leq&\frac{2(f-f(p)+1)(\langle\Phi\nabla\rho,\nabla\rho\rangle+\rho K)}{k(\rho^2+2)\log(\rho^2+2)}+\frac{4(f-f(p)+1)\rho^2\langle\Phi\nabla\rho,\nabla\rho\rangle}{k^2(\rho^2+2)^2[\log(\rho^2+2)]^2}.
\end{eqnarray*}
Now, since $|\nabla\rho|=1$ and $\Phi$ is positive semi-definite, one has
$\langle\Phi\nabla\rho,\nabla\rho\rangle\leq{\rm tr}(\Phi)$, so that the desired estimate
follows.
\end{proof}

\begin{corollary}\label{corol:functions bounded from above}
Let $\Sigma$ be a complete Riemannian manifold with sectional curvature
$K_{\Sigma}\geq 0$, and $f\in\mathcal D(\Sigma)$ be a function bounded from above. If $\Phi$ is positive
semi-definite and ${\rm tr}(\Phi)$ is bounded from above on $\Sigma$, then there exists
a sequence $(p_k)_{k\geq 1}$ in $\Sigma$ such that
\begin{equation}\label{eq:Omori-Yau for square 1}
f(p_k)>\sup_Mf-\frac{1}{k},\ \ |\nabla f(p_k)|<\frac{1}{k},\ \ \square f(p_k)<\frac{1}{k}.
\end{equation}
\end{corollary}

\begin{proof} Letting $C_1=\sup_{\Sigma}f$, it follows from (\ref{eq:bounds for the gradient of f})
that
\begin{eqnarray*}
|\nabla f(p_k)|&\leq&\frac{2(C_1-f(p)+1)}{k}\cdot\frac{\rho(p_k)}{\rho(p_k)^2+2}\cdot\frac{1}{\log(\rho(p_k)^2+2)}\\
&\leq&\frac{2(C_1-f(p)+1)}{k}\cdot\frac{1}{2\sqrt 2}\cdot\frac{1}{\log 2},
\end{eqnarray*}
so that
\begin{equation}\label{eq:limits of gradient}
\lim_{k\rightarrow+\infty}|\nabla f(p_k)|=0.
\end{equation}

If $f$ attains its maximum at some point of $\Sigma$, there is nothing to do. Otherwise, since
$(\Sigma,d)$ is a metric space, the sequence $(p_k)_{k\geq 1}$ whose existence is assured by the
previous theorem is such that $\lim_{k\rightarrow+\infty}\rho(p_k)=+\infty$. Hence,
since $K_{\Sigma}\geq 0$, it follows from lemma~\ref{lemma:bounds for K with bounded sectional curvature}
that, for sufficiently large $k$, one has $K(p_k)\leq\rho(p_k){\rm tr}(\Phi_{p_k})$.
Therefore, (\ref{eq:bounds for square of f}) gives
\begin{eqnarray*}
\square f(p_k)&\leq&\frac{2{\rm tr}(\Phi_{p_k})(C_1-f(p)+1)}{k}\left(\frac{\rho(p_k)^2+1}{\rho(p_k)^2+2}\right)\frac{1}{\log(\rho(p_k)^2+2)}\\
&&+\frac{4{\rm tr}(\Phi_{p_k})(C_1-f(p)+1)}{k^2}\left(\frac{\rho(p_k)}{\rho(p_k)^2+2}\right)^2\frac{1}{[\log(\rho(p_k)^2+2)]^2}.\\
&\leq&\frac{2C_2(C_1-f(p)+1)}{k\log 2}+\frac{C_2(C_1-f(p)+1)}{2k^2\log^22},
\end{eqnarray*}
so that
\begin{equation}\label{eq:limits of square of f}
\lim_{k\rightarrow+\infty}\square f(p_k)=0.
\end{equation}

The statement of the corollary follows from (\ref{eq:bounds for the gradient of f}),
(\ref{eq:limits of gradient}) and (\ref{eq:limits of square of f}), passing to a subsequence,
if necessary.
\end{proof}

\begin{corollary}\label{corol:functions bounded from below}
Let $\Sigma$ be a complete Riemannian manifold with sectional curvature
$K_{\Sigma}\geq 0$, and $f\in\mathcal D(\Sigma)$ be a function bounded from below. If $\Phi$ is positive
semi-definite and ${\rm tr}(\Phi)$ is bounded from above on $\Sigma$, then there exists
a sequence $(p_k)_{k\geq 1}$ in $\Sigma$ such that
\begin{equation}\label{eq:Omori-Yau for square 2}
f(p_k)<\inf_{\Sigma}f+\frac{1}{k},\ \ |\nabla f(p_k)|<\frac{1}{k},\ \ \square f(p_k)>-\frac{1}{k}.
\end{equation}
\end{corollary}

\begin{proof} Apply the previous corollary to $-f$.
\end{proof}

\section{Applications}\label{sec:applications}

Throughout this section, $\psi:\Sigma^n\rightarrow\overline M^{n+1}$ denotes, as before, a spacelike immersion
into a Lorentz manifold $\overline M$. In all that follows we set $\Phi=H_{r-1}P_{r-1}$, where $H_{r-1}$ and
$P_{r-1}$ are as in section~\ref{sec:preliminaries}. If $H_r=0$ on $\Sigma$, or else
$H_r>0$ on $\Sigma$ and $\psi$ has an elliptic point, then propositions~\ref{prop:ellipiticy_for_null_Sr}
and~\ref{prop:ellipiticy_for_positive_Sr} assure the semi-definiteness of $P_{r-1}$ (actually, $P_{r-1}$ is definite
when $H_r>0$).
Moreover, since
\begin{equation}
{\rm tr}\,\Phi=b_{r-1}H_{r-1}^2\geq 0,
\end{equation}
$\Phi$ is {\em positive} semi-definite in each of the above cases. In addition, if $H_{r-1}$ is
bounded on $\Sigma$, then the same is true of ${\rm tr}\,\Phi$, so that we can apply
corollaries~\ref{corol:functions bounded from above} and~\ref{corol:functions bounded from below}
to such a $\Phi$.

The following proposition is the analogue, in our context, of a lemma due to K. Akutagawa
(cf.~\cite{Akutagawa:87}).

\begin{proposition}\label{prop:Akutagawa for square operator}
Let $\overline M^{n+1}$ be a Lorentz manifold and $\psi:\Sigma^n\rightarrow\overline M^{n+1}$ a spacelike
immersion from a complete Riemannian manifold $\Sigma$ of sectional curvature $K_{\Sigma}\geq 0$ into $\overline M$.
Suppose that, for some $0<r\leq n$, $H_{r-1}$ is bounded on $\Sigma$ and one of the following is true:
\begin{enumerate}
\item[$(a)$] $H_r=0$ on $\Sigma$.
\item[$(b)$] $H_r>0$ on $\Sigma$ and $\psi$ has an elliptic point.
\end{enumerate}
If $f\in\mathcal D(M)$ is nonnegative and such that $\square f\geq af^{\beta}$, for some $a>0,\beta>1$, then $f\equiv 0$.
\end{proposition}

\begin{proof} Let $\phi:\mathbb R_+^*\rightarrow\mathbb R_+^*$ be a smooth function
to be chosen later, and $g=\phi\circ f$. Then $\nabla g=\phi'(f)\nabla f$ and
\begin{eqnarray*}
\square g&=&{\rm tr}(\Phi\,\text{Hess}\,g)=H_{r-1}L_{r-1}(g)=H_{r-1}\,\text{div}\,(P_{r-1}\nabla g)\\
&=&\phi'(f)H_{r-1}L_{r-1}(f)+\phi''(f)H_{r-1}\langle P_{r-1}\nabla f,\nabla f\rangle\\
&=&\phi'(f)\square f+\phi''(f)\langle\Phi\nabla f,\nabla f\rangle\\
&=&\phi'(f)\square f+\frac{\phi''(f)}{\phi'(f)^2}\langle\Phi\nabla g,\nabla g\rangle,
\end{eqnarray*}
so that
$$-\frac{\phi''(f)}{\phi'(f)^2}\langle\Phi\nabla g,\nabla g\rangle+\square g=\phi'(f)\square f.$$

Letting $\phi(t)=\frac{1}{(1+t)^{\alpha}}$, $\alpha>0$, one gets
$$\phi'(t)=-\alpha\phi(t)^{\frac{\alpha+1}{\alpha}},\ \ \frac{\phi''(f)}{\phi'(f)^2}=\left(\frac{\alpha+1}{\alpha}\right)\frac{1}{\phi(t)},$$
and hence
$$\left(\frac{\alpha+1}{\alpha}\right)\langle\Phi\nabla g,\nabla g\rangle-\phi(f)\square g=\alpha\phi(f)^{\frac{2\alpha+1}{\alpha}}\square f\geq a\alpha\frac{f^{\beta}}{(1+f)^{2\alpha+1}}.$$
If one now takes $\alpha=\frac{\beta-1}{2}>0$, we arrive at
\begin{equation}\label{eq:fundamental step}
\left(\frac{\alpha+1}{\alpha}\right)\langle\Phi\nabla g,\nabla g\rangle-g\square g\geq a\alpha\left(\frac{f}{1+f}\right)^{\beta}.
\end{equation}

Since $g$ is bounded from below, by corollary~\ref{corol:functions bounded from below}
we get a sequence $(p_k)$ of points in $M$ such that
$$g(p_k)<\inf_Mg+\frac{1}{k},\ \ |\nabla g|(p_k)<\frac{1}{k},\ \ \square g(p_k)>-\frac{1}{k}.$$
Therefore, $f(p_k)\rightarrow\sup_Mf$, and taking into account that
$$\langle\Phi\nabla g,\nabla g\rangle\leq({\rm tr}\,\Phi)|\nabla g|^2=b_{r-1}H_{r-1}^2|\nabla g|^2,$$
we get from (\ref{eq:fundamental step}) that
$$b_{r-1}H_{r-1}^2\left(\frac{\alpha+1}{\alpha k^2}\right)-\frac{1}{k}\left(\inf_Mg+\frac{1}{k}\right)\geq a\alpha\left(\frac{f(p_k)}{1+f(p_k)}\right)^{\beta}.$$

Making $k\rightarrow+\infty$, we get $\sup_Mf=0$, and since $f\geq 0$ this gives $f\equiv 0$.
\end{proof}

Let $M^n$ be a connected, $n$-dimensional oriented Riemannian manifold and $I\subset\mathbb R$ an interval.
In the product manifold $\overline M^{n+1}=I\times M^n$, let $\pi_I$ and $\pi_M$ denote the projections onto the
$I$ and $M$ factors, respectively. If $g:I\rightarrow\mathbb R$ is a positive smooth function, we obtain a
particular class of Lorentz metrics in $\overline M^{n+1}$ by setting
$$\langle v,w\rangle_p=-\langle(\pi_I)_*v,(\pi_I)_*w\rangle+(f\circ\pi_I)(p)^2\langle(\pi_M)_*v,(\pi_M)_*w\rangle,$$
for all $p\in\overline M$ and all $v,w\in T_p\overline M$. Furnished with such a metric, $\overline M$ is called
a Generalized Robertson-Walker (GRW) spacetime, and will be denoted by writing $\overline M^{n+1}=-I\times_gM^n$.
In Cosmology, a GRW gives a simple, physically plausible relativistic model (cf.~\cite{O'Neill:83}),
so a natural space to work with.

In a GRW spacetime $\overline M^{n+1}=-I\times_gM^n$ one has the globally defined conformal vector field
$V=g\partial_t$, which is even {\em closed}, in the sense that its dual $1-$form is closed; moreover, one
can easily prove that ${\rm div}\,V=(n+1)g'$. If $\psi:\Sigma^n\rightarrow\overline M^{n+1}$ is a spacelike
immersion, we oriented $\Sigma$ by choosing a timelike unit normal vector field $N$.
For future use, we quote lemma 5.4 of~\cite{AliasBrasil:03}, where the reader can also find a thorough
discussion of a class of spacetimes more general than that of GRW's.

\begin{lemma}\label{lemma:existence of elliptic point}
Let $\overline M^{n+1}=-I\times_gM^n$ be a GRW spacetime, and $\psi:\Sigma^n\rightarrow\overline M^{n+1}$ a
spacelike immersion. If the restriction of $g\circ\pi_I$ to $\psi(\Sigma)$ attains a local minimum at some
$p\in\psi(\Sigma)$, such that $g'(\pi_I(p))\neq 0$, then $p$ is an elliptic point for $\Sigma$.
\end{lemma}

The following proposition is due to L.J. Al\'{\i}as and A.G. Colares, as lemma $4.1$ of
preprint~\cite{AliasColares:06}. Here, and for the sake of completeness, we present a more direct proof.

\begin{proposition}\label{prop:Lr of height function}
Let $\overline M^{n+1}=-I\times_gM^n$ be a GRW spacetime, and $\psi:\Sigma^n\rightarrow\overline M^{n+1}$ a
spacelike immersion. If $h=\pi_{I_{|\Sigma}}:\Sigma^n\rightarrow I$ is the height function of $\Sigma$, then
\begin{equation}\label{eq:Lrh}
L_r(h)=-(\log f)'\{b_rH_r+\langle P_r(\nabla h),\nabla h\rangle\}-b_rH_{r+1}\langle N,\partial_t\rangle.
\end{equation}
\end{proposition}

\begin{proof} One has
\begin{eqnarray*}
\nabla h&=&\nabla(\pi_{I_{|\Sigma}})=(\overline\nabla\pi_I)^{\top}=-\partial_t^{\top}\\
&=&-\partial_t-\langle N,\partial_t\rangle N,
\end{eqnarray*}
where $\overline\nabla$ denotes the gradient with respect to the metric of the ambient space and
$X^{\top}$ the tangential component of a vector field $X\in\mathcal X(\overline M)$ in $\Sigma$.
Now fix $p\in M$, $v\in T_pM$ and let $A$ denote the Weingarten map with respect to $N$. Write
$v=w-\langle v,\partial_t\rangle\partial_t$, so that $w\in T_p\overline M$ is tangent to the
fiber of $\overline M$ passing through $p$. By repeated use of the formulas of item (2) of proposition $7.35$
of~\cite{O'Neill:83}, we get
\begin{eqnarray*}
\overline\nabla_v\partial_t&=&\overline\nabla_w\partial_t-\langle v,\partial_t\rangle\overline\nabla_{\partial_t}\partial_t=\overline\nabla_w\partial_t\\
&=&(\log f)'w=(\log f)'(v+\langle v,\partial_t\rangle\partial_t).
\end{eqnarray*}
Thus,
\begin{eqnarray*}
\nabla_v\nabla h&=&\overline\nabla_v\nabla h+\langle Av,\nabla h\rangle N\\
&=&\overline\nabla_v(-\partial_t-\langle N,\partial_t\rangle N)+\langle Av,\nabla h\rangle N\\
&=&-(\log f)'w-v(\langle N,\partial_t\rangle)N+\langle N,\partial_t\rangle Av+\langle Av,\nabla h\rangle N\\
&=&-(\log f)'w+(\langle Av,\partial_t\rangle-\langle N,\overline\nabla_v\partial_t\rangle)N+\langle N,\partial_t\rangle Av+\langle Av,\nabla h\rangle N\\
&=&-(\log f)'w+(\langle Av,\partial_t^{\top}\rangle-\langle N,(\log f)'w\rangle)N+\langle N,\partial_t\rangle Av+\langle Av,\nabla h\rangle N\\
&=&-(\log f)'w-(\log f)'\langle v,\partial_t\rangle\langle N,\partial_t\rangle N+\langle N,\partial_t\rangle Av\\
&=&-(\log f)'\{v-\langle v,\partial_t\rangle(-\partial_t-\langle N,\partial_t\rangle N)\}+\langle N,\partial_t\rangle Av\\
&=&(\log f)'(-v+\langle v,\partial_t^{\top}\rangle\nabla h)+\langle N,\partial_t\rangle Av\\
&=&-(\log f)'(v+\langle v,\nabla h\rangle\nabla h)+\langle N,\partial_t\rangle Av.
\end{eqnarray*}
Now, by fixing $p\in\Sigma$ and an orthonormal frame $\{e_i\}$ at $T_p\Sigma$, one gets
\begin{eqnarray*}
L_rh&=&{\rm tr}({\rm Hess}\,h)=\sum_{i=1}^n\langle\nabla_{e_i}\nabla h,P_{r}e_i\rangle\\
&=&\sum_{i=1}^n\langle-(\log f)'(e_i+\langle e_i,\nabla h\rangle\nabla h)+\langle N,\partial_t\rangle Ae_i,P_{r}e_i\rangle\\
&=&-(\log f)'\{{\rm tr}(P_r)+\langle P_r(\nabla h),\nabla h\rangle\}+\langle N,\partial_t\rangle{\rm tr}(AP_r).
\end{eqnarray*}
The result follows from the formulas for the traces of $P_r$ and $AP_r$.
\end{proof}

Now we consider a particular model of Lorentzian GRW, the {\em Steady State space}, namely
\begin{equation}\label{eq:steady state space}
\mathcal H^{n+1}=-\mathbb R\times_{e^t}\mathbb R^n.
\end{equation}
This spacetime corresponds to the steady state model of the universe proposed by Bondi, Gold and Hoyle
(cf.~\cite{Hawking:73}, chapter $5$).

A spacelike immersion $\psi:\Sigma^n\rightarrow\mathcal H^{n+1}$ such that $H_r=0$ on $\Sigma$ is said to be
$r-$maximal. If $h\geq t_0$ on $\Sigma$ for some $t_0\in\mathbb R$, $\psi$ is said to be a spacelike hypersurface
{\em over the slice} $M_{t_0}=\{t_0\}\times M$.

\begin{theorem}\label{thm:application 1}
There exists no $r$-maximal complete spacelike hypersurface $\psi:\Sigma^n\rightarrow\mathcal H^{n+1}$ over the slice
$M_{t_0}$ of $\mathcal H^{n+1}$, with sectional curvature $K_{\Sigma}\geq 0$ and such that
$C_1\leq |H_{r-1}|\leq C_2$, for some positive constants $C_1,C_2$.
\end{theorem}

\begin{proof}
Suppose, by contradiction, the existence of such a hypersurface. Given a smooth function
$\varphi:\mathbb R\rightarrow\mathbb R$, a straightforward computation shows that
$$L_r(\varphi\circ h)=\varphi^{\prime\prime}(h)\left\langle P_{r}\nabla h,\nabla h\right\rangle + \varphi^{\prime}(h)L_{r}(h),$$
so that equation (\ref{eq:Lrh}) gives
$$L_{r-1}\left(e^{-h+t_0}\right)=e^{-h+t_0}\left\{\langle P_{r-1}\nabla h,\nabla h\rangle+b_{r-1}H_{r-1}\right\}.$$
Consequently, since $\Phi=H_{r-1}P_{r-1}$, we have that
\begin{eqnarray*}
\Box(e^{-h+t_0})&=&{\rm tr}(\Phi\,{\rm Hess}(e^{-h+t_0}))=H_{r-1}L_{r-1}(e^{-h+t_0})\\
&=&e^{-h+t_0}\left\{2\langle\Phi(\nabla h),\nabla h\rangle+b_{r-1}H_{r-1}^{2}\right\}.
\end{eqnarray*}
Thus, since $\Phi$ is positive semi-definite and $h-t_0\geq0$, we get
\begin{equation*}
\Box(e^{-h+t_0})\geq C_1^2b_{r-1}e^{\beta(-h+t_0)}, \ \forall\,\beta>1.
\end{equation*}
Therefore, from proposition~\ref{prop:Akutagawa for square operator} we conclude that
$e^{-h+t_0}\equiv 0$, which is an absurd.
\end{proof}

\begin{remark}
As a consequence of Bonnet-Myers theorem, a complete spacelike hypersurface $\psi:\Sigma^n\rightarrow\mathcal H^{n+1}$
having (not necessarily constant) mean curvature $H$ satisfying $|H|\leq\rho<2\sqrt{n-1}/n$
($\rho$ a real constant), has to be compact; in fact, for such a bound on $H$, Gauss' equation would give
$${\rm Ric}_{M}\geq(n-1)-n^2\rho^2/4>0,$$
where ${\rm Ric}_{\Sigma}$ denotes the Ricci curvature of $\Sigma$. However, since the Steady State space is not
spatially closed, i.e., since its Riemannian fiber is not compact, such a hypersurface does not exist
(cf. proposition $3.2(i)$ of~\cite{AliasRomeroSanchez:95}).
\end{remark}

\begin{remark}
As a special case of the reasoning of the above remark, we see that there are no complete maximal spacelike
hypersurfaces in $\mathcal H^{n+1}$. Theorem~\ref{thm:application 1} can thus be seen as a sort of generalization of
this situation for higher order mean curvatures.
\end{remark}

In what follows, we say that a spacelike hypersurface $\psi:\Sigma^n\rightarrow\mathcal H^{n+1}$ has the
{\em same time-orientation} of $\partial_t$ if $\Sigma$ is oriented by the choice of a timelike unit normal vector field
$N$, such that $\langle N,\partial_t\rangle\leq-1$; otherwise we say that $\Sigma$ has time-orientation {\em opposite}
to that of $\partial_t$.

\begin{theorem}\label{thm:application 2}
Let $\psi:\Sigma^n\rightarrow\mathcal H^{n+1}$ be a complete spacelike hypersurface
over a slice $M_{t_0}$ of $\mathcal H^{n+1}$, with sectional curvature $K_{\Sigma}\geq 0$ and time-orientation
opposite to that of $\partial_t$. If $H_r>0$ and
$C_1\leq H_{r-1}\leq C_2$ for some positive constants $C_1$ and $C_2$,
then the height function $h=\pi_{\mathbb{R}_{|\Sigma}}$ does not attain a local minimum on $\Sigma$.
\end{theorem}

\begin{proof}
Suppose that, for some such hypersurface $\psi:\Sigma^n\rightarrow\mathcal H^{n+1}$, the height function do attains
a local minimum, at $p\in\psi(\Sigma)$, say. Since $g=e^h$ on $\psi(\Sigma)$, $p$ is also a
local minimum for $g\circ\pi_I$, and hence lemma~\ref{lemma:existence of elliptic point} assures the existence
of an elliptic point for $\psi(\Sigma)$; therefore, by proposition~\ref{prop:ellipiticy_for_positive_Sr} $P_{r-1}$
is positive definite.

Now, equation (\ref{eq:Lrh}) gives
$$L_{r-1}\left(e^{-h+t_0}\right)=e^{-h+t_0}\left\{\langle P_{r-1}\nabla h,\nabla h\rangle+b_{r-1}[H_{r-1}+H_{r}\langle N,\partial_t\rangle]\right\}.$$
Thus, taking once more $\Phi=H_{r-1}P_{r-1}$ we get
\begin{eqnarray*}
\Box(e^{-h+t_0}) &=&{\rm tr}(\Phi\,{\rm Hess}(e^{-h+t_0}))=H_{r-1}L_{r-1}(e^{-h+t_0})\\
&=&e^{-h+t_0}\{2\langle\Phi(\nabla h),\nabla h\rangle+b_{r-1}[H_{r-1}^2+H_{r-1}H_r\langle N,\partial_t\rangle]\}.
\end{eqnarray*}
Since $\Phi$ is positive definite, $\langle N,\partial_t\rangle\geq1$ and $h-t_0\geq0$, we finally obtain
\begin{equation*}
\Box(e^{-h+t_0})\geq C_1^2b_{r-1}e^{\beta(-h+t_0)},\ \forall\,\beta>1.
\end{equation*}
Therefore, by proposition~\ref{prop:Akutagawa for square operator} we conclude that
$e^{-h+t_0}\equiv0$, which is an absurd.
\end{proof}

When $r=2$ and $\Sigma$ has time-orientation opposite to that of $\partial_t$, lemma 3.2
of~\cite{AliasColares:06} assures the ellipticity of $L_1$ whenever $H_2>0$. Since, by Gauss' equation, this is the
same as asking that $\Sigma$ has scalar curvature $R<n(n-1)$, one can reason as in the previous result to obtain the
following

\begin{theorem}\label{thm:application 3}
There exists no complete spacelike hypersurface $\psi:\Sigma^n\rightarrow\mathcal H^{n+1}$ over a slice
$M_{t_0}$ of $\mathcal H^{n+1}$, with sectional curvature $K_{\Sigma}\geq 0$ and satisfying the following conditions:
\begin{enumerate}
\item[(a)] $\Sigma$ has scalar curvature $R<n(n-1)$;
\item[(b)] If the time-orientation of $\Sigma$ is opposite to that of $\partial_t$, then its mean curvature $H$
is such that $C_1\leq H\leq C_2$, for some positive constants $C_1$ and $C_2$.
\end{enumerate}
\end{theorem}

\bibliographystyle{amsplain}

\end{document}